\documentclass[12pt]{article}         
\usepackage{amsfonts,amssymb}
\makeatletter
\def\bdraft{\pagestyle{myheadings} 
           \textheight=10.5truein \textwidth=7.5truein \parindent=8pt
           \voffset=-1truein \topmargin=30pt \headheight=10pt \headsep=3pt
           \ifcase \@ptsize \hoffset=-1.5truein \or \hoffset=-1.35truein
                        \or \hoffset=-1.15truein \fi}
\def\quality{\textheight=240mm \textwidth=160mm \topmargin=0Truein
             \ifcase \@ptsize \hoffset=-23mm
                     \or \hoffset=-20mm \or \hoffset=-15mm \fi}
\makeatother 
\quality 

\def\beq#1#2{\begin{equation} \label{#1} #2 \end{equation}}
\def\bea#1{\begin{eqnarray*} #1 \end{eqnarray*}} \def\a{\!\!\!&\!\!\!\!&}
\def\beaq#1#2{\label{#1} \begin{eqnarray} #2 \end{eqnarray}}
\def\n{\noindent}   \def\map{T}   
    
\def\IR{{\mathbb{R}}}  \def\IZ{{\mathbb{Z}}}  
\def\toas#1{\stackrel{#1}{\longrightarrow}}
  \def\phi{\varphi}   \def\la{\lambda}
  \def\mod1{\,({\rm mod\ } 1)\,}
\def\t{\tilde} \def\rho{\varrho} 

\def\function#1{\left\{\!\!\!\begin{array}{ll} #1 \end{array} \right.}

\def\proof{\smallskip \noindent {\bf Proof. \ }}       
\def\blanksquare{\,\,\,$\sqcup\!\!\!\!\sqcap$}         
\def\qed{\hfill\blanksquare\linebreak\smallskip\par}   

\def\thname{Theorem}     \def\lmname{Lemma}      \def\prname{Proposition}
\def\dfname{Definition}  \def\crname{Corollary}  \def\rmname{Remark}

\newtheorem{theorem}{\thname}
\newtheorem{lemma}{\lmname}
\newtheorem{corollary}[lemma]{\crname}   

\newtheorem{dftn}{\dfname}[section]
\newtheorem{rmrk}[lemma]{\rmname}

             \catcode`@=11 \@addtoreset{equation}{section} \catcode`@=12

\makeatletter\def\fps@figure{htbp}\makeatother 
\newcommand\mlbscale{1pt} 
\newif\iffigs\figstrue 
\def\bline(#1,#2)(#3,#4)(#5){\put(#1,#2){\line(#3,#4){#5}}}  

\def\bfig(#1,#2)#3#4{\begin{figure} \begin{center}
    \framebox{\setlength{\unitlength}{\mlbscale}
       \iffigs \begin{picture}(#1,#2) #3 \end{picture}
       \else \begin{picture}(60,10)(0,0)
                   \put(0,0){\framebox(60,10){Figure}} \end{picture} \fi}
    \end{center} \caption{#4} \end{figure}}

\def\Bfig(#1,#2)#3#4{\begin{figure} \begin{center}
    \setlength{\unitlength}{\mlbscale}
       \iffigs \begin{picture}(#1,#2) #3 \end{picture}
       \else \begin{picture}(60,10)(0,0)
                   \put(0,0){\framebox(60,10){Figure}} \end{picture} \fi
    \end{center} \caption{#4} \end{figure}}

\def\bpic(#1,#2)#3{\setlength{\unitlength}{\mlbscale}
    \begin{picture}(#1,#2) #3 \end{picture}}

\def\bB{{\bf B}}   
 \def\cB{{\cal B}} \def\cC{{\cal C}}
\def\bR{{\bf R}} \def\rho{\varrho} \def\cM{{\cal M}}

\def\*#1{#1^*}    \def\0#1{\breve#1}  \def\2#1{\acute#1}
\def\G{\Delta}    \def\b#1{\overline{#1}}   \def\bp{\bar p}
\def\intp#1{\left\lfloor#1\right\rfloor}

\def\?#1{} 

\begin{document}
\title{Stochastic stability of traffic maps}
\author{Michael Blank\thanks{
        Russian Academy of Sci., Inst. for
        Information Transmission Problems, ~
        e-mail: blank@iitp.ru}
        \thanks{This research has been partially supported
                by Russian Foundation for Basic Research.}
       }
\date{October 8, 2012} 
\maketitle

\begin{abstract}%
We study ergodic properties of a family of traffic maps acting in
the space of bi-infinite sequences of real numbers. 
The corresponding dynamics mimics the motion of vehicles in 
a simple traffic flow, which explains the name. Using connections to 
topological Markov chains we obtain nontrivial invariant measures, 
prove their stochastic stability, and calculate the topological entropy. 
Technically these results in the deterministic setting are related to 
the construction of measures of maximal entropy via measures 
uniformly distributed on periodic points of a given period, while 
in the random setting we directly construct (spatially) Markov 
invariant measures. In distinction to conventional results the 
limiting measures in non-lattice case are non-ergodic. 
Average velocity of individual ``vehicles'' as a function of their 
density and its stochastic stability is studied as well.
\end{abstract}%

\bigskip

\n 2000 Mathematics Subject Classification. 
Primary: 28D05; Secondary:  28D20, 37A50 34F05 35B35.

\n Key words: dynamical system, topological Markov chain, 
invariant measure, stochastic stability, traffic flow.

\section{Introduction}\label{s:intro}

Speaking about stochastic stability of a dynamical system one means
(see e.g. \cite{Si,Bl-mon}) 
that the most important statistical quantities related to the
dynamics (e.g. Sinai-Bowen-Ruelle invariant measures) depend
continuously on the addition of a small amount of true random noise.
We study ergodic properties of a family of traffic maps acting in
the space of bi-infinite sequences. The corresponding dynamics
mimics the motion of vehicles in a simple traffic flow, which
explains the name. Our aim is to construct nontrivial invariant
measures of the traffic maps and to show that they are stable with
respect to ``natural'' random perturbations -- by which we mean that
the motion of individual vehicles are performed with a probability
$p$. Then the case $p=1$ corresponds to the deterministic traffic
map, while situations with $p<1$ may be considered as random
perturbations.

We introduce the following notation. Under an
{\em admissible configuration} $x^t$ at time $t\in\IZ_+\cup\{0\}$
we mean an ordered countable set of particles (balls) of radius
$r\ge0$, centers of which are located at the points
$x^t:=(x_i^t)_{i\in\IZ}\subset\bR\subseteq\IR^1$ such that
$$x_i^t + r\le x_{i +1}^t-r.$$
The set of all admissible configurations we denote by $X=X(r,\bR)$.
By $v>0$ denote the maximal possible movement of a single particle
per unit time, i.e. 
$$0\le x_i^{t+1}-x_i^t \le v.$$ 
The parameter $p\in(0,1]$ stands for the probability of movement for 
individual particles. A single particle performs a totally asymmetric 
random walk (which explains why the processes of this sort are often 
called TASEP - Totally Asymmetric Exclusion Process) with jumps of 
size $v$, occurring with probability $p$, until its motion does not 
interfere with the motions of other particles.

Exclusion processes (collective random walks of countable
collections of particles with hard core interactions) introduced by
Frank Spitzer in 1970 appear naturally in a broad list of scientific
applications starting from various models of traffic flows
\cite{NS,GG,ERS,Bl-erg,Bl-hys}, molecular motors and protein
synthesis in biology, surface growth or percolation processes in
physics (see \cite{Pe,BFS} for a review), and up to the analysis of
Young diagrams in Representation Theory \cite{CMOS}. Continuous time
versions of these processes are reasonably well understood (see e.g.
\cite{Lig} for a general account and \cite{An,AAV,EFM,FM} for recent
results). The main difficulty in the analysis of discrete time
versions of exclusion processes is that an arbitrary (infinite)
number of particle interactions may happen simultaneously. To
overcome this difficulty one needs to develop principally new
approaches. Note that in the one-dimensional setting under
consideration the basic restriction is the preservation of the order
of particles, i.e. the particles can not overtake each other under
dynamics.

We consider discrete time Markov processes $\pi(p,v,r,\bR)$ acting
in the space of configurations $X=X(r,\bR)$ and the dynamics of
individual particles in the configuration is defined by the relation
\beq{e:dyn}{x_i^{t+1}=\function{\min\{x_i^t+v,x_{i+1}^t-2r\}
                          &\mbox{with probability } p \\
                     x_i^t &\mbox{with probability } 1-p} .}

\def\particle#1{
     \put(0,0){\circle{30}} \bline(0,0)(0,-1)(25) \put(-3,-32){$x_{#1}^t$}
     \bline(0,20)(1,0)(15) \put(5,22){$r$}
     \put(0,35){\vector(1,0){70}} \put(20,39){$v_{}$}}
\Bfig(150,50)
      {\footnotesize{
       \thicklines
       \bline(0,10)(1,0)(150)
       \put(30,10){\particle{i}} \put(110,10){\particle{i+1}}
       \thinlines \bline(46,28)(1,0)(48) \put(65,30){$\Delta_i$}
       \bezier{50}(30,10)(55,-20)(78,10) \put(54,-5){\vector(1,0){4}} \put(53,-1){$p$}
       \put(78,10){\circle{30}} \bline(78,10)(0,-1)(25) \put(75,-22){$x_i^{t+1}$}
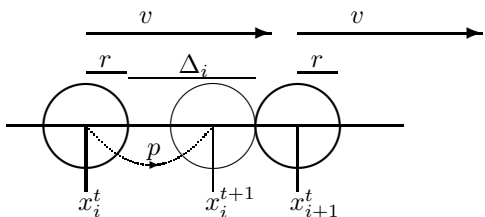
      } } {Exclusion process in continuum \label{f:tasep-c}}

Consider three, at first glance, very different types of exclusion
processes satisfying the relation (\ref{e:dyn}). 
Processes of type 1 act on a lattice $\bR=\IZ^1, r=1/2, v\in\IZ_+^1$. 
One lattice site may be occupied by at most one
particle. Models of this type are widely used to describe the motion
of vehicles on a single-lane road (see for example \cite{NS,SS}).

Processes of type 2 also act on a lattice
$\bR=\IZ^1,~ v\in\IZ_+^1$, but $r=0$. The fundamental difference
is that a lattice site can be occupied by an arbitrary number of
particles. Models of this sort with a continuous-time are called
zero-range processes (see, for example. \cite{EH05}) and it is
convenient to use them to simulate a communication line, in which
particles represent equal packets of information, waiting in queues
to communication servers, located at sites of the lattice $\bR$.
In terms of quantum statistical mechanics processes of type 1 and
2 are related as interacting Fermi gas and free Bose gas.

Processes of type 3, which are a special case of exclusion type
processes introduced in \cite{Bl10}, act not on a lattice but on the
continuous space $\bR=\IR^1, r\ge0, v\in\IR_+^1$.
It will be shown (Lemma~\ref{l:proc-relations}) that for $r=1/2$
these processes contain all realizations of the processes of the 1st
type, and for $r=0$ all realizations of the processes the 2nd type.
This allows one to obtain simultaneously analytical results for all
cases.
The consideration of exclusion type processes in continuum not only
simplifies the analysis but due to the presence of additional
symmetries (absent in lattice versions) it offers a possibility to
obtain new results about the lattice cases unavailable otherwise.

A {\em density}~ of a configuration $x^t$ (the number of particles 
per unit length) is defined as 
$$\rho(x^t):=\lim_{n\to\infty}n/(x_{n-1}^t-x_0^t),$$ 
if the latter limit makes sense.\footnote {The one-sidedness of the
   growth of segments $(x_0^t,x_{n-1}^t)$ is due to the fact that
   all particles move in the same direction. The definition of the
   density in the general case is more complicated
   (see e.g. \cite {Bl10}).} %
As we shall see the density is the first integral of the processes 
under study and thus the phase space is foliated naturally by invariant 
subsets $X^{(\rho)}$ consisting of admissible configurations of 
density $\rho$. 

The formula~(\ref{e:dyn}) describes the processes through dynamics 
of configurations of ordered particles. This is convenient for the 
analysis of the particles motion, but does not allow the study of 
invariant measures (stationary distributions). To this end we 
consider a modification of the process under study in which the 
particles are indistinguishable from each other. Speaking about 
invariant measures we always refer to the invariant measures of 
this modification. The set of probability $\pi(p,v,r,\bR)$-invariant 
measures we denote by $\cM_{p,v,r,\bR}$.

To date, the mathematical description of invariant measures for the 
deterministic processes under study is essentially absent, and the only 
classification result \cite{BF} gives just a formal description of invariant 
measures for the simplest lattice process $\pi(p=1,v=1,r\in\{0,1/2\},\IZ)$ 
and says nothing even about the existence of nontrivial (non-atomic) 
invariant measures. Note that the processes $\pi(p=1,v,r,\bR)$ possess 
periodic trajectories of all possible periods and hence atomic measures 
supported by these trajectories. 
In Theorems~\ref{t:mes-main},\ref{t:nonatom-det},\ref{t:lattice-det} 
for each density $\rho$ we prove the existence of a non-atomic invariant 
measure $\mu_\rho$ supported by configurations of the given density. 
Moreover we give an explicit construction in terms of (spatially) Markov 
measures related so subshifts of finite type. 

Having a large number of invariant measures it is important to 
distinguish ``physically relevant'' ones. In the case of low dimensional 
dynamical systems one often uses for this purpose the concept of 
Sinai-Bowen-Ruelle (SBR) measures. Roughly speaking the latter means 
that for a reasonably large family of ``good'' initial probabilistic measures 
(say all absolutely continuous measures) their images under dynamics 
converge in Cesaro means to the same SBR measure. In the present 
infinite-dimensional setting the choice of ``good'' initial measures is 
not obvious and the control over the convergence of their images is not 
available at present. An alternative approach consists in the 
analysis of stochastic stability of invariant measures. In a number of 
cases it has been shown (see, e.g., \cite{Bl-mon}) that SBR measures 
are exactly those that are stochastically stable. 

As we already noted a natural choice of random perturbations for the 
deterministic processes $\pi(p=1,v,r,\bR)$ are random processes 
$\pi(p<1,v,r,\bR)$. Here again in the discrete time case only very partial 
results about invariant measures are known in the literature. 
Nevertheless the simplest lattice case $\pi(p<1,v=1,r=1/2,\IZ)$ and 
especially its much simpler continuous time version was intensively 
studied in physics literature from this point of view 
(see \cite{An,EFM,EH05,FM,Gr,NS,SS,SSNI}). In particular, in \cite{SS,SSNI} 
using a mean field approximation and an interesting and nontrivial 
combinatorial argument (apparently not quite complete without the exact 
analysis of the limit construction) the authors have found necessary conditions 
(equivalent to our formula~(\ref{e:id-markov}) for the existence of nontrivial 
invariant measures. On the other hand, the next step -- the proof that 
the constructed measures are indeed invariant under dynamics  was not 
explicitly addressed. 
The construction proposed in the proof of Lemma~\ref{l:markov} 
provides not only a complete proof of this statement but also in its 
first few lines gives the crucial formula~(\ref{e:id-markov}) without 
any limit transitions. By means of a completely different approach 
M.Kanai \cite{Ka} has studied the same process on a finite discrete
ring (instead of the infinite lattice) and obtained an exact rather
complicated formula for the invariant measure as a function of the
ring's length $L$. When $L$ goes to infinity this invariant measure 
converges to the limit derived in \cite{SS}, which is of independent 
interest despite a number of restrictions assumed in this paper. 
The approaches elaborated in \cite{SS,SSNI,Ka} use in one way 
or another properties of invariant measures of the process
$\pi(p,v=1,r=1/2,\{1,2,\dots,L\})$ on a finite ring. Actually in
\cite{SS,SSNI} the Markov structure of the invariant measure was
discovered (without the proof of its existence) while in \cite{Ka}
the unique invariant measure is no longer Markovian and becomes
Markovian only in the limit as $L\to\infty$. A coupling construction 
proposed by L.~Gray \cite{Gr} allows one to show that for a fixed density
 $\rho$ the Markov invariant measure is ergodic and  unique among
translationally invariant measures.

It is worth noting that the above mentioned lattice constructions cannot 
be extended to the case of long jumps ($v>1$). To overcome this 
difficulty we first study the continuous setting to obtain the invariant 
measure for the process $\pi(p<1,v=1,r=1/2,\IR)$ 
and then extend this result for all values of $v,r$ using 
a wide set of symmetries available in $\IR$ (but not in $\IZ$). 
After that we are able to return back to the lattice cases.

Despite the fact that the constructions of invariant measures in 
deterministic and stochastic settings are rather different and neither 
of them can be applied in another setting, the proof of stochastic 
stability is relatively simple and is based on the explicit description 
of the invariant measures. 

Our main result about invariant measures is as follows.

\begin{theorem}\label{t:mes-main} The process $\pi(p,v,r,\IR)$ for 
each $0<p\le1, v>0,r\ge0$ and each particle density $\rho\in(0,1/(2r)]$ 
possesses a nontrivial invariant measure $\mu_\rho^p$ supported 
by the set of configurations of density $\rho$ and 
$\mu_\rho^p\toas{p\to1}\mu_\rho^1$ in the weak sense.
\end{theorem}

One of the most important characteristics of the processes under
consideration is an {\em average velocity of particles} during time $t>0$:
$$V(x,i,t):=(x_i^t-x_i^0)/t.$$ In \cite{Bl10} it was shown that
under very general assumptions (including, in particular, dynamics
of particles moving in opposite directions), which definitely hold
true for all processes considered in the present work, the a.s.
limit of the statistics $V(x,i,t)$ as $t\to\infty$ (if it exists)
depends only on the density $\rho$ of the configuration $x$
(see Theorem~\ref{t:vel-part}). Therefore, to calculate the average
velocity $$V(x)=V(\rho(x)):=\lim_{t\to\infty}V(x,i,t)$$ (where the
convergence is considered in the almost sure sense) for a configuration
of a given density  it is enough to determine the latter for a specially
selected initial configuration of the same density (for which the
calculation is convenient). Note that the existence of the limit under
question, not to mention the explicit formulas for it, previously has
been proven only in the deterministic setting (i.e. for $p=1$), see \cite{Bl10}.
The following result gives explicit relations between the average velocity 
and other parameters of the process $\pi(p,v,r,\IR)$. 

\begin{theorem}\label{t:vel-main} For the process %
$\pi(p,v,r,\IR)~~\forall v,\rho\in\IR_+^1, ~0\le r<1/(2\rho),~p\in(0,1)$ %
with any admissible initial configuration of density $\rho$
the average velocity $V(\rho,r)$ is well defined and is
calculated as follows:~  %
\beq{e:vel-change}{ V(\rho,r)=(1-2r\rho)~V(\frac{\rho}{1-2r\rho},0),} %
\beq{e:vel}
{V(\rho,0)=\left[1+v\rho-\sqrt{(1+v\rho)^2-4pv\rho}\right]/(2\rho)
\toas{p\to1}\min\{1/\rho,v\}. }%
\end{theorem}

\Bfig(120,100)
      {\small{
       \thicklines
       \bline(0,0)(1,0)(140) \bline(0,0)(0,1)(90)
       \bline(0,73)(1,0)(60)
       \bezier{200}(60,73)(75,5)(135,2)
       \thinlines
       \bezier{30}(60,73)(60,36)(60,0)  \put(30,78){$p=1$}
       \put(145,-5){$\rho(x)$}  \put(-24,85){$V(x)$}
       \put(-7,70){$v$}
       \put(85,30){$1/\rho(x)$} \put(-8,-5){$0$}
       \put(55,-12){$1/v$}
       \bezier{100}(0,65)(55,65)(60,50)
       \bezier{200}(60,50)(75,5)(135,1) \put(30,45){$p<1$}
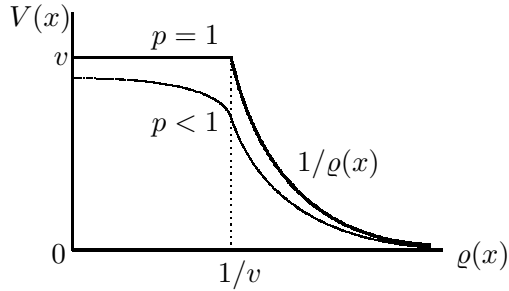
      }      }{Fundamental diagrams (average velocity against density)
                 for the process $\pi(p,v,r=0,\IR)$.  \label{f:tasep-c-f}}

In physics literature results similar to Theorem~\ref{t:vel-main}
are known for the simplest lattice TASEP process
$\pi(p<1,v=1,r=1/2,\IZ)$ (see \cite{SS,SSNI,Ka}). Here again 
one needs to make comments about missing arguments of the same 
sort as in the analysis of invariant measures plus that the result should 
hold not only for ``typical'' initial configurations (see discussion below 
and in Section~\ref{s:vel}). 

The dynamical coupling construction developed in \cite{Bl10,Bl11}
allows one to get complete information about the properties of average
velocities in the deterministic setting (i.e. for $p=1$).
This construction does not require the study of (numerous) invariant
measures of the process, but gives only conditional (upon their existence)
results in the stochastic case (albeit under much broader assumptions
about the process: the local velocities $v=v_i$ are iid random variables).
In the stochastic case one cannot avoid the analysis of invariant
measures. Explicit expression for the invariant measure $\mu$ allows 
one to derive (by Birkhoff's ergodic theorem) the formula for the average 
velocity for $\mu$-a.a. initial configurations. To extend this result to 
a much more broad setting (formulated in Thereom~\ref{t:vel-main})   
for all particle configurations for which densities are well defined, one 
needs to use metric properties of the process obtained earlier. Namely 
in \cite{Bl10} it has been shown that for all particle configurations of 
a given density the average velocity is the same. Therefore it is enough 
to calculate this statistics for a single suitable initial configuration (e.g. 
a configuration typical with respect to the measure $\mu$).

The paper is organized as follows. In Section~\ref{s:metric} we
discuss connections between the processes under study and briefly
review known results related to particle densities and average velocities.
Ergodic properties of the deterministic traffic map will be discussed 
in Section~\ref{s:erg-det}, which we shall finish by the calculation 
of the topological entropy of the traffic map in continuum. The true 
random setting (i.e. $0<p<1$) will be analyzed in Section~\ref{s:erg-stoch}. 
In Section~\ref{s:vel} we shall study average particle velocities,  
obtain explicit relations between these statistics and parameters 
of the process, and prove their stochastic stability. 
Finally in Section~\ref{s:het-size} we shall discuss heterogeneous 
versions of the processes under consideration (particles with different 
sizes in the same configuration and presence of randomly
distributed static obstacles in space).

\section{Metric properties}\label{s:metric}

Here we shall discuss connections between the processes under study
and briefly review known results related to particle densities and
average velocities.

Let $\G_i=\G_i(x^t):=x_{i+1}^t-x_i^t-2r$  stands for the distance
between boundaries of the balls corresponding to subsequent particles
in the configuration $x^t$(see Fig.~\ref{f:tasep-c}). We refer to
$\G_i$ as the {\em gap} corresponding to the $i$-th particle
in $x^t$. If $\G_i(x^t,z)<v$ we say that the $i$-th particle is {\em
blocked} (meaning that its motion is blocked by the $i+1$-th
particle) at time $t$ and {\em free} otherwise. By a {\em cluster}
of particles in a configuration $x^t\in X$ we mean a locally maximal
collection of consecutive blocked particles.

To emphasize the dependence of various statistics on the parameter
$r\ge0$ (ball's radius) we denote by $x(r)$ and $X(r,\bR)$ the configurations
with balls of radius $r$ and the corresponding space of admissible
configurations. In the special case of point-particles ($r=0$) the
dependence on the radius will be omitted.

We say that two particle processes whose dynamics is described by
the relation~(\ref{e:dyn}) are {\em statically coupled} if all
random choices related to the motion of particles with the same
indices in these processes coincide.

The following results show relations between the three types
of exclusion processes (two lattice and one in continuum) introduced 
in the previous Section.

\begin{lemma}\label{l:proc-relations}
For any given $v,p$ the following relations are valid: %
\begin{enumerate}
\item $\pi(p,v\in\IZ_+,r=1/2,\IZ)=\pi(p,v,r=1/2,\IR)$ if $x^0(r=1/2)\subset\IZ$; %
\item $\pi(p,v\in\IZ_+,r=0,\IZ)=\pi(p,v,r=0,\IR)$ if $x^0(r=0)\subset\IZ$; %
\item $\forall r>0$ there exists an affine homeomorphism
      $\phi=\phi_r:X(r,\IR)\to X(0,\IR)$, such that %
      $$\phi\circ\pi(p,v,r,\IR)=\pi(p,v,r=0,\IR)\circ\phi;$$ %
\item $\forall ~u,v>0,~0<p\le1$ a sub-lattice
      $\bR_{u,v}:=v\IZ+u$ is invariant with respect to the process
      $\pi(p,v,r=0,\IR)$,  i.e.
      $x^0\subset\bR_{u,v}~\Rightarrow~x^t\subset\bR_{u,v}
       ~\forall t>0$.
\end{enumerate}  \end{lemma}

\proof The only nontrivial statement here is the item~(3).
Observe that any two particle configurations $x(r),~\2x(\2r)$
having the same sequence of gaps $\G:=\{\G_i\}$ may be transformed
to each other by a one-to-one map %
\beq{e:r->r'}{\2x_i(\2r)=\phi_r(x_i(r)):=x_i(r)-2i(r-\2r)~~\forall i\in\IZ.} %
This allows one to choose the homeomorphism $\phi_r$ as follows
$$(\phi_r(x(r)))_i:=x_i(r)-2ir~~\forall i\in\IZ.$$
Now the conjugation between the processes having particles of
different sizes follows from the equality between sequences of gaps
$\Delta^t$ between particles in the statically coupled processes
$\pi(p,v,r,\IR)$ and $\pi(p,v,0,\IR)$. \qed

The correspondence between densities and average particle velocities
for statically coupled processes with configurations consisting of balls
of different radiuses $r\ne\2r\ge0$ is summarized as follows.

\begin{lemma} \cite{Bl10} \label{l:den-vel-cor}
\begin{enumerate}
\item The density $\rho(x^t)$ is preserved by dynamics, i.e.
     $\rho(x^t)=\rho(x^{t+1})~~\forall t\ge0$.
\item Let configurations $x(r)\in X(r,\IR),~r>0$ and
    $\2x(\2r)\in X(\2r,\IR),~\2r\ge0$ have the same sequence of gaps
    $\{\Delta_i\}$. Then
    $$\rho(x(r))=\frac{\rho(\2x(\2r))}{1+2(r-\2r)\rho(\2x(\2r))} .$$
\item Let additionally the processes $x^t(r)$ and $\2x^t(\2r)$ be
    statically coupled. Then $\forall i,t$ $$V(x(r),i,t)=V(\2x(\2r),i,t).$$
\end{enumerate}
\end{lemma}%

Results obtained in \cite{Bl10} allow also to claim the coincidence
of average velocities both for individual particles and for configurations
having equal densities.

\begin{theorem}\label{t:vel-part} (\cite{Bl10}) For the process
$\pi(p,v,r,\bR)~~\forall v,r\in\bR, ~\rho>0, ~p\in(0,1]$
the following claims are valid. Let $x,y\in X(r,\bR)$,
$\rho(x)=\rho(y)=\rho$ and let the average velocity $V(y)$ be well
defined. Then $$|V(x,i,t)-V(y,j,t)|\toas{t\to0}0$$ for any $i,j\in\IZ$.
\end{theorem}

\section{Ergodic properties: deterministic case (p=1)}\label{s:erg-det}

Despite the abundance of invariant measures supported by
time-periodic trajectories (existing for all periods) nothing was
known previously about the existence of nontrivial invariant measures.
Our main result in this direction in the case $\bR:=\IR$ is the
following claim.

\begin{theorem}\label{t:nonatom-det}
The deterministic process $\pi(p=1,v,r,\IR)$ for each $v>0,r\ge0$ and
each particle density $\rho\in(0,1/(2r)]$ possesses a non-atomic
invariant measure $\mu_\rho$ such that $E_{\mu_\rho}[x^{t}_i]=\rho$.
\end{theorem}

The proof of this result will be given step by step in 
Sections~\ref{s:lat-traf},\ref{s:massiv-det}. 
In the simplest lattice case the construction of the invariant measure
will be described explicitly in terms of the so called (spatially) Markov 
measures. The notion {\em Markov measure} here refers to the unique
invariant measure for the Markov shift on $\{0,1\}^\IZ$ with a 
non-degenerate transition matrix
$(p_{ij}),~i,j\in\{0,1\},~p_{ij}\ge0,~\sum_j p_{ij}=1$.
In terms of configurations of particles, 1 corresponds to the presence
of a particle at a site of the lattice $\IZ$, and 0 - its absence.

It is worthwhile giving an equivalent definition of the Markov measure
without any connection to random processes (see e.g. \cite{To}).
A measure $\mu$ on  $\{0,1\}^\IZ$ is Markov iff for any cylinder\footnote{a
   set of all sequences taking a given collection of values at a given
   collection of indices, called its base}  %
with the base $[AbC]$, where $A$ and $B$ two arbitrary finite
binary sequences while $b$ is a single binary letter, we have 
\beq{e:markov}{\mu([b])\cdot\mu([AbC]) = \mu([Ab])\cdot\mu([bC]) .}
Note that \cite{To} claims that the Markov measures appear extremely 
rarely in cellular automata systems and thus their existence for traffic maps 
turns out to be a surprise both in deterministic and random settings.

\subsection{Lattice traffic map}\label{s:lat-traf}

\begin{theorem}\label{t:lattice-det}
The deterministic lattice process $\pi(p=1,v=1,r=1/2,\IZ)$ possesses
a Markov invariant measure $\mu$, being a weighted sum
of two measures of maximal entropy for shift-maps acting in
opposite directions.
Besides, there is a 1-parameter family of invariant measures
$\{\mu_{1,1,1/2,\IZ,\rho}\}_\rho$ such that
$E_{\mu_{1,1,1/2,\IZ,\rho}}[x^{t}_i]=\rho$.
\end{theorem}

Note that the construction of Markov invariant measures which we shall
elaborate for the analysis of the random setting in Section~\ref{s:erg-stoch}
cannot be applied here and we shall use a very different strategy related
to the construction of measures of maximal entropy well known in 
theory of low-dimensional hyperbolic dynamical systems 
(see e.g. \cite{Bowen}).

Let $\bB:=\{0,1\}^\IZ$ be the space of binary sequences and let
$\bB_\pm$ be its two subspaces defined by the transition matrices
$M_+:=\left(\begin{array}{cc}1&1\\1&0 \end{array}\right)$ and
$M_-:=\left(\begin{array}{cc}0&1\\1&1 \end{array}\right)$,
namely 
$$ \bB_+:=\{b\in\bB:~~ b_i+b_{i+1}\ne2~\forall i\in\IZ \} ,$$
$$ \bB_-:=\{b\in\bB:~~ b_i+b_{i+1}\ne0~\forall i\in\IZ \} .$$
Consider two shift-maps $\sigma_\pm:\bB\to\bB$ defined as follows:
$$(\sigma_+b)_i:=b_{i-1}, \qquad (\sigma_-b)_i:=b_{i+1} ,$$
i.e. these two maps shift a sequence by one position in opposite
directions. Obviously 
$$\sigma_\pm\bB_\pm\equiv\bB_\pm.$$

Denote by $\map:\bB\to\bB$ the map corresponding to the
deterministic process $\pi(p=1,v=1,r=1/2,\IZ)$.
Recall that the density of a binary sequence $b\in\bB$ is defined
exactly as the density of a configuration of ones in $b$ given in
Section~\ref{s:intro}.

Theorem~\ref{t:lattice-det} claims only the existence of Markov
invariant measures but in fact we shall prove that the restriction
of the dynamical system $(\map,\bB)$ to the Cantor sets $\bB_\pm$
possesses {\em massive} invariant measures in the sense that these
measures are positive on each open subset. Before to give the proof
of this result let us discuss reasons for the absence of massive
invariant measures in the entire space. In Section~\ref{s:erg-stoch} 
we shall show that in the true random case ($0<p<1$) for each 
$a=p_{01}\in(0,1)$ there is a nontrivial measure $\mu_a^{(p)}$ 
(constructed in the proof of Theorem~\ref{t:mes-mas-1}) whose 
value on the cylinder with the base $[1100]$ is equal to
$$ \mu_a^{(1)}([1100]) = p_1p_{11}p_{10}p_{00}
     =  p_1(1-p_{10})p_{10}p_{00} .$$
On the other hand, by (\ref{e:p10})
$$ p_{10}=(1-p_{01})/(1-pp_{01}) \toas{p\to1}1 .$$
Therefore for each $p_{01}$ we have
$$ \mu_a^{(1)}([1100])\toas{p\to1} 0 ,$$
and thus the limiting (as $p\to1$) measure cannot be massive. 
This explains why in the deterministic case when $p=1$ one should
not expect the existence of true massive invariant measures. We
conjecture that for each $\map$-invariant probability measure $\eta$
there exists an open subset $B\in\bB$ with $\eta(B)=0$. However
at present we do not have a complete proof of this statement.

\bigskip

\n{\bf Proof of Theorem~\ref{t:lattice-det}}.
Let us start with a recipe of the construction of measures of
maximal entropy for the shift-map $\sigma:\bB_M\to\bB_M$ with an
irreducible transition matrix $M=(m_{ij})_{i,j\in\{0,1\}}$
(see e.g. \cite{Bowen}).
For each $n\in\IZ_+$ denote by $\mu_n$ the probability measure
uniformly distributed on points of period $n$ of the map $\sigma$.
Then the measures $\mu_n$ weakly converge as $n\to\infty$ to
a massive (on $\bB_M$) probability measure $\mu_\sigma$.
The latter coincides with the unique invariant distribution for the
stationary Markov chain with the transition probability matrix
$$P_\sigma:=(m_{ij}m_j/(\la_M m_i)),$$
where $\la_M$ is the maximal eigenvalue of the matrix $M$ and
$(m_i)$ is the corresponding eigenvector. The measure $\mu_\sigma$
represents the only invariant measure maximizing the metric
entropy of the dynamical system $(\sigma,\bB_M)$ with the transition
matrix $M$, which explains the reason for its name. Calculating the
leading elements of the spectra of the matrices $M_\pm$ we get
$$\la_+=\la_-=:\la:=(1+\sqrt5)/2, \quad m_0^+=m_1^-=1/\la,
\quad m_1^+=m_0^-=1-1/\la.$$

Despite the absence of a similar result for the dynamical system
$(\map,\bB)$ we try to follow this recipe. A serious problem here
is that we do not have good control over all time-periodic points of
a given period for the map $\map$. To overcome this difficulty we
select a subset of time-periodic points with which we shall work.

Observe that if a configuration $b\in\bB$ is spatially periodic then
$\forall t\in\IZ_+~\map^tb$ is spatially periodic with the same period.
This fact follows immediately from the definition of the map $\map$
and is discussed in detail in \cite{Bl-erg}. Additionally in \cite{Bl-erg}
it has been shown that if either $\sigma_+^nb=b$ and $b\in\bB_+$
or $\sigma_-^nb=b$ and $b\in\bB_-$ then this spatially periodic point
is time-periodic for the map $\map$ with the same period $n$.
Note that
$\rho(b)\le1/2$ if $b\in\bB_+$ and $\rho(b)\ge1/2$ if $b\in\bB_-$.

The absence of non spatially periodic time-periodic points of the
map $\map$ would complete the description of time-periodic points.
Unfortunately this is not the case. The point is that the map is not
one-to-one. Using this let us construct a sketch of a counter-example.
Let $b$ be spatially and time periodic for the map $\map$ with the
same period $n$, i.e. $b_i=b_{i+n}~\forall i\in\IZ$ and $\map^nb=b$.
Thus a finite word $(b_{n},b_{-n+1},\dots,b_{-1})$ is a pre-image
of the word $(b_0,b_1,\dots,b_{n-1})$ under the action of the map
$\map^n$. Since the map $\map$ is not bijective there exists another
pre-image $(b'_{-n},b'_{-n+1},\dots,b'_{-1})$ of the word
$(b_0,b_1,\dots,b_{n-1})$. Similarly we consider a pre-image of the
word $(b'_{n},b'_{-n+1},\dots,b'_{-1})$  under the action of the map
$\map^n$ and denote it by $(b'_{-2n},b'_{-2n+1},\dots,b'_{-n-1})$.
Continuing this procedure we are getting a point
$b':=(\dots,b'_{-2},b'_{-1},b_0,b_1,b_2,\dots)\in\bB$ such that
$b'$ is no longer spatially periodic but $\map^nb'=b'$.

Nevertheless we can use a spatially periodic part of time-periodic points
for our construction. For each $n\in\IZ_+$ consider a probability measure
$\mu_n^*$ uniformly distributed on spatially and time periodic points
of period $n$ of the map $\map$. As we already noted each spatially
and time periodic point of the map $\map$ is spatially periodic point of
one of the shift-maps $\sigma_\pm$. Due to the symmetry of the
motion of ones and zeros under the action of the map $\map$ the
number $C_n^+$ of $n$-periodic points of density less or equal to
$1/2$ differs at most by $n$  from the number $C_n^-$ of
$n$-periodic points of density larger than $1/2$, while $C_n^\pm$
are of order $n^\la$ (see the calculation of the exponent  $\la>0$ below). 
Therefore one can represent the measure $\mu_n^*$ as follows
$$\mu_n^*=(C_n^+\mu_n^+ + C_n^-\mu_n^-)/(C_n^+ + C_n^-) ,$$
where $\mu_n^\pm$ are probability measures uniformly distributed
on $n$-periodic points of the shift-maps $\sigma_\pm$. Being
uniformly distributed on time-periodic trajectories the measure $\mu_n^*$
is $\map$-invariant for each $n$.

We already know that the measures $\mu_n^\pm$ converge
as $n\to\infty$ to the corresponding measures of maximal entropy
$\mu_\sigma^\pm$ for the shift-maps $\sigma_\pm$. Therefore
$$ \mu_n^*\toas{n\to\infty}(\mu_\sigma^+ + \mu_\sigma^-)/2=:\mu_\map .$$
The weak massive property for the limit measure follows from the similar
statement for the measures $\mu_\sigma^\pm$ on the sets $\bB_\pm$.
Indeed by the construction for each cylinder on $\bB_+$ its
$\mu_\sigma^+$ measure is positive, and the similar statement holds
for each cylinder on $\bB_-$ and the measure $\mu_\sigma^-$.

Note that according to our construction all three limit measures
$\mu_\map,\mu_\sigma^\pm$ are $\map$-invariant. Moreover,
the equality of the leading eigenvalues $\la_+=\la_-=:\la:=(1+\sqrt5)/2$
implies the coincidence of the corresponding metric entropies
(being equal to $\ln\la$) .

It remains to construct for each $\rho\in[0,1]$ the $\map$-invariant
measure $\mu_\rho$ supported by the configurations of density
$\rho$. This can be done using either Gibbsian reconstruction of the
measures $\mu_\sigma^\pm$ or by an explicit representation in
terms of Markov shifts. Let us discuss the latter approach.

The action of a shift-map $\sigma$ is equivalent to the time-shift along
realizations of the Markov chain with the transition probability matrix
compatible with the corresponding binary transition matrix $M$. Recall that
matrices $A=(a_{ij})$ and $B=(b_{ij})$ with nonnegative entries are
{\em compatible} if $a_{ij}b_{ij}=0$ implies $a_{ij}+b_{ij}=0$.

A probability transition matrix compatible with the matrix $M_+$ is
written as
$P_+:=\left(\begin{array}{cc}1-a&a\\1&0 \end{array}\right)$
with a single parameter $0\le a\le1$. The function $\rho:=a/(1+a)$
defines a bijection between the values of the parameter $a$ and 
the set of particle densities $\rho$. Thus for each density
$\rho\in[0,1/2]$ we obtain a massive (spatially) Markov
$\map$-invariant measure $\mu_\rho$.

Similarly one considers  the matrix $M_-$ which allows one to construct
massive invariant measures $\mu_\rho$ with $\rho\in(1/2,1]$.
\qed

Theorem~\ref{t:lattice-det} gives a complete recipe for the
construction of nontrivial invariant measures for the process
$\pi(p=1,v=1,r=1/2,\IZ)$.

\subsection{Existence of massive invariant measures in continuum}\label{s:massiv-det}

In this Section we develop a special machinery extending the measures
$\mu_\rho$ (constructed for the lattice case) first to nontrivial invariant 
measures of the deterministic process $\pi(p=1,v=1,r=1/2,\IR)$ acting
on the real line and then to deterministic processes 
$\pi(p=1,v>0,r\ge0,\IR)$ with arbitrary $v>0,r\ge0$.

Recall that in this more complicated `continuous' setting a finite
{\em cylinder} with the base defined by a finite subset of integers
$I$ and a collection $C:=\{C_i\}_{i\in I}$ of open
intervals\footnote{In general the cylinder $\cC_{I,C}$
   might be empty for nonempty sets $I,C$.} %
is the subset $\cC_{I,C}:=\{x\in X:~~x_i\in C_i~~\forall i\in I\}$.
We endow the space of admissible configurations $X$ with the
$\sigma$-algebra $\cB$ generated by the finite cylinders defining a
topology in this space.

By Lemma~\ref{l:proc-relations}(1) the measures
$\mu_\rho=\mu_{\rho,r}$
are $\pi(p=1,v=1,r=1/2,\IR)$-invariant. For each $r'\ge0$ applying
the affine transformation (\ref{e:r->r'}), obtained in
Lemma~\ref{l:proc-relations}(3), to the measure $\mu_\rho$
we are getting a new probability measure $\mu_{\rho',r'}$ with
$\rho':=\rho/(1+2(1/2-r')\rho)$, which is
$\pi(p=1,v=1,r',\IR)$-invariant. Since $\rho$ takes all
values from the interval $[0,1]$ the new variable $\rho'$ takes all
values from the interval $[0,1/(2r')]$. 

Making yet another spatial change of variables $z\to vz+w$ with
parameters $v>0,~w\ge0$ we are obtaining from the measure
$\mu_{\rho',r'}$ a two-parameter family of probability measures
$\mu(p=1,\rho'',r'',v,w)$ supported on configurations of balls of
radius $r'':=vr'$ and having density $\rho'':=\rho'/v$. Applying again
Lemma~\ref{l:proc-relations}(4) we see that for each $\rho'',r',v,w$
the measure $\mu(p=1,\rho'',r'',v,w)$ is $\pi(p=1,v,r'',\IR)$-invariant.

On the other hand, these measures cannot be massive, since they
are supported by very thin sets. To construct a massive
(albeit non-ergodic) invariant measure from this family we consider
a measure $\mu(p=1,\rho'',r'',v)$ having marginals $\mu(p=1,\rho'',r'',v,w)$
on sub-lattices $v\IZ+w$ and uniformly distributed with respect to the
parameter $w\in[0,v)$.

This finishes the proof of the existence of massive invariant measures 
for the deterministic traffic map in continuum.

The construction of the massive invariant measure in the lattice case
($\bR:=\IZ$) with general $v\in\IZ_+$ is similar except that
$w\in\{0,1,\dots,v-1\}$ and $r"=1/2$ or $r"=0$ to be able to work
with point particles in lattice setting. 

\subsection{Entropy}\label{s:entropy}

To finalize the description of ergodic properties of the processes
under study in the pure deterministic setting we show (following the
approach developed in \cite{Bl10}) that these processes are
strongly chaotic. Our dynamical system is defined by a deterministic
map $\map_v:X\to X$ from the set of admissible configurations into itself.
Our aim is to show that the topological entropy of this map is
infinite.\footnote{Normally one says that a map is chaotic if
   its topological entropy is positive, so infinite value of
   the entropy indicates a very high level of chaoticity.} %

We refer the reader to \cite{Bi,Wa} for detailed definitions of
the topological and metric entropies for deterministic dynamical
systems and their properties that we use here. To avoid
difficulties related to the non-compactness of the phase space we
define the topological entropy of a map $\map_v$ (notation
$h_{{\rm top}}(\map_v)$) as the supremum of metric entropies of
this map taken over all probabilistic invariant measures (compare
to the conventional definition of the topological entropy and its
properties discussed, e.g. in \cite{Wa}).

We start the analysis with the action of a shift-map in continuum
$\sigma_v:X\to X$ defined as %
$$(\sigma_v x)_i:=x_i+v~~~ i\in\IZ, ~x\in X,~v>0.$$

\begin{lemma}\label{l:entropy-shift-map}
The topological entropy of the shift-map in continuum $\sigma_v$
is infinite.
\end{lemma}
\proof The continuity of the shift-map in continuum in the topology 
induced by the $\sigma$-algebra $\cB$ generated by finite cylinders 
is implied by the fact that a preimage of a finite cylinder under 
the action of $\sigma_v$ is again a finite cylinder. 

The idea of the calculation of the entropy is to construct an invariant 
subset of $X$ on which the map $\sigma_v$ is isomorphic to the full 
shift-map in the space of sequences with a countable alphabet. 
The result follows from the observation that the topological entropy of the
full shift-map $\sigma^{(n)}$ with the alphabet consisting of $n$
elements is equal to $\ln n$ (see, e.g. \cite{Bi,Wa}).

Let $\alpha:=\{\alpha_i\}_{i\in\IZ_+}$ with $\alpha_i\in(0,v)$ and
let $\alpha^{(n)}:=\{\alpha_i\}_{i=1}^n$. Consider a sequence of
subsets $X^{(n)}\subset X$ consisting of {\em all} configurations
$x\in X$ satisfying the condition 
$$x_{2k}\in{v}\IZ, ~x_{2k+1}\in x_{2k} + \alpha^{(n)}~~\forall k\in\IZ.$$ 
Then $X^{(n)}$ is $\sigma_v$-invariant and the restriction $\sigma_v|X^{(n)}$ is
isomorphic to the full shift-map $\sigma^{(n)}$ with the alphabet
$A^n$ consisting of $n$ elements $\{a_i\}$ of type %
$a_i:=\{[0,\alpha_i),[\alpha_i,v)\}$, i.e. each element is
represented by a pair of neighboring intervals. Therefore the
topological entropy of $\sigma^{(n)}$ is equal to %
$\ln n\toas{n\to\infty}\infty$. \qed

\begin{theorem}\label{t:entropy-exclusion}
The topological entropy of the traffic map in continuum $\IR$ is infinite.
\end{theorem}
\proof The traffic map is continuous in the topology induced by the 
$\sigma$-algebra $\cB$ generated by finite cylinders by the same 
argument as in the case of the shift-map.

Observe that the subset %
$$X_{>v}:=\{x\in X:~~\G_i(x)\ge v~~ \forall i\in\IZ\}$$ 
of the set of admissible configurations is $\map_v$-invariant.
Therefore 
$$h_{{\rm top}}(\map_v)\ge h_{{\rm top}}(\map_v|X_{>v})$$ 
and for our purposes it is enough to show that the latter is
infinite. On the other hand, by the definition of the map
$\map_v$ we have $\map_v|X_{>v}\equiv\sigma_v|X_{>v}$.

One cannot apply the result of Lemma~\ref{l:entropy-shift-map} 
directly because in the case under consideration the gaps 
between particles are greater or equal to $v$ by the construction, 
while in the proof of Lemma~\ref{l:entropy-shift-map} the gaps 
did not exceed $v$. To this end one sets $\alpha_i\in(v,2v)$ and
modifies the definition of $X^{(n)}$ as follows: 
$$\tilde{X}^{(n)}:=\{x_{2k}\in3v\IZ, \quad 
                              x_{2k+1}\in x_{2k} + \alpha^{(n)}\quad 
                              \forall k\in\IZ \} .$$ %
Consider the the alphabet $A^{(n)}$ with elements of type
$a_i:=\{[0,\alpha_i),[\alpha_i,3v)\},~~\alpha_i\in\alpha^{(n)}$. 
Then the $3$-d power of the map $\map_v|X_{>v}$ ~ is 
isomorphic to the full shift-map $\sigma^{(n)}$
with the alphabet  $A^{(n)}$. Using that 
$$3h_{{\rm top}}(\map_v|X_{>v})=h_{{\rm top}}((\map_v|X_{>v})^3)
=h_{{\rm top}}(\sigma^{(n)})=\ln n$$ we get the result. \qed

\section{Ergodic properties: stochastic case ($0<p<1$)}\label{s:erg-stoch}

\subsection{The simplest lattice process $\pi(p,v=1,r=1/2,\IZ)$}
\label{s:erg-simple-stoch}

We start with the construction of a nontrivial Markov invariant measure for 
the simplest lattice TASEP.

\begin{lemma}\label{l:markov}
Let $P:=\left(\begin{array}{cc}
                  p_{00}&p_{01}\\p_{10}&p_{11} \end{array}\right)$
be a probability matrix with positive entries and the left leading
normalized eigenvector $(p_0,p_1)$. Then the measure $\mu$
on $X$, defined on cylinders by the relation
\beq{e:markov-mes}{\mu([a_1,a_2,\dots,a_n])
       := p_{a_1}\prod_{i=1}^{n-1}p_{a_i,a_{i+1}},~~a_i\in\{0,1\}, }
is invariant with respect to the process $\pi(p,v=1,r=1/2,\IZ)$ iff
\beq{e:id-markov}{p_{00}p_{11}=(1-p) p_{10}p_{01} .} 
\end{lemma}

\proof For a Markov chain $\pi$ acting on a space of binary sequences
we say that a cylinder $B$ is a (partial) {\em pre-image} of a cylinder $A$
if the probability $\Pr(\pi(B)=A)>0$. We need to check that for any
cylinder its measure is equal to the sum of measures of all its pre-images
multiplied by the corresponding transition probabilities. The proof
follows by induction on the cylinder's length $n$.

We start with cylinders of length 1 and 2. The following tables show all
pre-images of cylinders of length 1 and 2. The left column corresponds
to cylinders and their pre-images, while the right column shows stationary
probabilities for the cylinder (the 1st line) and stationary probabilities
for its pre-images multiplied by the transition probabilities.
Here $\bp:=1-p$ is the probability that the jump does not take place.

\bigskip

\def\btab#1{\begin{tabular}{l|l} #1 \end{tabular}}

\btab{~0 & $p_0$ \\ \hline
          00 & $p_0p_{00}$ \\
          10 & $\bp p_1p_{10}$\\
         ~10& $p p_1p_{10}$\\}
\quad
\btab{~~1 &  $p_1$ \\ \hline
          ~10 &  $\bp p_1p_{10}$\\
          ~11 &   $p_1p_{11}$\\
          10  &   $pp_1p{10}$\\}
\qquad \hskip-0.5mm
\btab{~10  & $p_1p_{10}$ \\ \hline
         ~10  & $\bp p_1p_{10}$ \\
         ~110 & $pp_1p_{11}p_{10}$ \\
          100  & $pp_1p_{10}p_{00}$ \\
          1010 & $p^2p_1p_{10}^2p_{01}$ \\}
\quad
\btab{~11   & $p_1p_{11}$ \\ \hline
         ~110 & $\bp p_1p_{11}p_{10}$ \\
         ~111 & $p_1p_{11}^2$ \\
         1010  & $p\bp p_1p_{10}^2p_{01}$ \\
         1011  & $pp_1p_{10}p_{01}p_{11}$ \\}
\bigskip \par
\btab{~00   & $p_0p_{00}$ \\ \hline
         ~000 & $p_0p_{00}^2$ \\
         ~100 & $\bp p_1p_{10}p_{00}=\bp p_0p_{01}p_{00}$ \\
         0010 &  $pp_0p_{00}p_{01}p_{10}=pp_1p_{10}^2p_{00}$ \\
         1010 &  $\bp pp_1p_{10}^2p_{01}=\bp pp_0p_{01}^2p_{10}$ \\}
\quad
\btab{~01  & $p_0p_{01}=p_1p_{10}$ \\ \hline
         ~10 & $pp_1p_{10}$ \\
        0010 & $\bp p_0p_{00}p_{01}p_{10}=\bp p_1p_{10}^2p_{00}$ \\
        0011 & $p_0p_{00}p_{01}p_{11}=\bp p_1p_{10}^2p_{01}$ \\
        1010 & $\bp^2 p_1p_{10}^2p_{01}=\bp p_1p_{10}p_{00}p_{11}$ \\
        1011 & $\bp p_1p_{10}p_{01}p_{11}$ \\}

\bigskip

In the first two cases the equivalence of the probability in the 1st line to
the sum of other probabilities is trivial. In the 3d case assuming that
$pp_{10}\ne0$ we get the following condition for the equivalence  %
\beq{e:id2}{ 1=\bp+pp_{11}+pp_{00}+p^2p_{10}p_{01} %
    \Longrightarrow ~1=p_{00}+p_{11}+pp_{10}p_{01} ,}%
while in the 4th case assuming that $p_1p_{11}\ne0$  we have
$$ 1=\bp p_{10}+p_{11}+p\bp p_{10}^2p_{01}/p_{11}+pp_{10}p_{01} %
    \Longrightarrow ~p_{00}p_{11}=\bp p_{10}p_{01} .$$

The last calculation proves that the assumption~(\ref{e:id-markov}) is necessary.
A direct calculation shows that despite appearances the relation~(\ref{e:id2})
is equivalent to (\ref{e:id-markov}).

The checking of the cases 5 and 6 can be done similarly using the property
$p_0p_{01}=p_1p_{10}$ (which we already used in the corresponding tables).

Assume now that the claim is already proven for all cylinders with bases less
or equal to $n>1$. To reduce the case of length $n+1$ to $n$ observe that
the only difference between smaller cylinder's length to the larger one consists
of an additional letter $0$ or $1$ at the right end of the cylinder's base.
Therefore it is enough to show that the corresponding probability changes by
$p_{ab}$, where $a$ is the last letter in the shorter cylinder and $b$ is the
new letter. To prove this we consider all 4 possibilities:
$$1\to10, 1\to11, 0\to00, 0\to01.$$

The 2nd and 3d situations are relatively simple:

$1\to11: \qquad  10\to110, 11\to111$ and $10\to1010+1011$
$$ pp_1p_{10}(\bp p_{10}p_{01}+p_{01}p_{11})
  =pp_1p_{10}(p_{00}p_{11}+p_{01}p_{11})=pp_1p_{10}p_{11}.$$

$0\to00: \qquad 00\to000, 10\to100$ and $.10\to0010+1010$\par
(here and in the sequel ``$.$'' stands for an arbitrary symbol)
$$pp_1p_{10}(p_{10}p_{00}+\bp p_{01}p_{10})
=pp_1p_{10}(p_{10}p_{00}+p_{00}p_{11})
=pp_1p_{10}p_{00}.$$

Note that we use that the transition probabilities depend only
on the previous letter.

In the 1st situation we proceed as follows ~
$$1\to10: \quad 10+11\to10+110, \quad 10\to100+1010.$$
Denoting the product of the measure of a cylinder $A$ and the corresponding
transition probability by $K(A)$ we get
\bea{
K([10])+K([11])\a=p_1(\bp p_{10}+p_{11})=p(p_{10}-pp_{10}+p_{11})
                =p_1(1-pp_{10}) \\
K([10])+K([110])\a=p_1p_{10}(\bp+pp_{11})=p_1p_{10}(1-p+pp_{11}) \\
                \a=p_1p_{10}(1-p(1-p_{11}))=p_1p_{10}(1-pp_{10}).}
Thus the 2nd sum differs from the 1st by the desired multiplier $p_{10}$.

The 4th situation is a bit more complicated. In this case we split the
3d pre-image $.10$ of 0 into two parts 010 and 110 giving contributions
$$K([010])=pp_0p_{01}p_{10}, \quad K([110])=pp_1p_{11}p_{10}$$
respectively (observe that
$pp_0p_{01}p_{10}+pp_1p_{11}p_{10}=pp_1p_{10}$),
and the 1st pre-image 10 of $01$ also splits into two parts 010 and 110
giving contributions 
$$K([010])=pp_0p_{01}p_{10}=pp_1p_{10}^2, \quad
   K([110])=pp_1p_{11}p_{10}$$ 
respectively. Then we gather them as follows: ~ %
$$00+010\to010+0010+0011, \quad 10+110\to110+1010+1011.$$
The following simple calculation checks the correctness of this
construction:
\bea{
K([00])+K([010])\a=p_0(p_{00}+pp_{01}p_{10})=p_0(1-p_{11})=p_0p_{10}  \\
K([0010])+K([0011])+K([010])\a=p_1p_{10}^2(\bp p_{00}+\bp p_{01}+p) \\
\a=p_1p_{10}^2=p_0p_{01}p_{10}=(p_0p_{10})p_{01} .}
\qed

\begin{theorem}\label{t:mes-mas-1}  The process
$\pi(p,v=1,r=1/2,\IZ)~~\forall v\in\IZ_+^1, p,\rho\in(0,1)$
possesses a 1-parameter family of probabilistic 
invariant measures
$\{\mu_{p,1/2,1,\IZ,\rho}\}_\rho$ positive on each open set 
and supported by the set of configurations of density $\rho$.
\end{theorem}

\proof 
Observe that the probability measure $\mu$ constructed in
Lemma~\ref{l:markov} due to the additional condition~(\ref{e:id-markov})
can be parametrized by a single parameter $a:=p_{01}\in(0,1)$.
From (\ref{e:id-markov}) we obtain
\beq{e:p10}{ p_{10}=(1-a)/(1-pa) .} %
In view of the equality between the density of particles $\rho$
and the stationary probability of ones $p_1$ under the Markov shift
with the transition matrix $P$ we have
$$\rho=a/(a+p_{10})=a(1-pa)/(1-pa^2). $$
Solving the last equality with respect to the parameter $a$, we get %
\beq{e:par}{a=\frac{1-\sqrt{1-4p\rho(1-\rho)}}{2p(1-\rho)}.} %
Therefore, for a given $p$ the constructed family of measures
$\{\mu_\rho^{(p)}\}_\rho$  is uniquely indexed by the density $\rho$.

The property that the measures $\{\mu_\rho^{(p)}\}_\rho$ are massive
follows immediately from the observation that the measure of an
arbitrary finite cylinder is positive. \qed

\subsection{Existence of massive invariant measures in continuum}\label{s:massiv}

Existence of massive $\pi(0<p<1,v=1,r=1/2,\IZ)$-invariant measures
$\mu_\rho^{(p)}$ for each $\rho$ is already proven in
Section~\ref{s:erg-simple-stoch}.

By Lemma~\ref{l:proc-relations}(1) the measures $\mu_\rho^{(p)}$ are
$\pi(p,v=1,r=1/2,\IR)$-invariant. For each $r'\ge0$ applying the
affine transformation (\ref{e:r->r'}), obtained in
Lemma~\ref{l:proc-relations}(3), to the measure $\mu_\rho^{(p)}$ we
are getting a new probability measure $\mu_{\rho',r'}$ with
$\rho':=\rho/(1+2(1/2-r')\rho)$, which is
$\pi(p,v=1,r',\IR)$-invariant. Since $\rho$ takes all values from
the interval $[0,1]$ the new variable $\rho'$ takes all values from
the interval $[0,1/(2r')]$.

Making yet another spatial change of variables $z\to vz+w$ with
parameters $v>0,~w\ge0$ we are obtaining from the measure
$\mu_{\rho',r'}$ a two-parameter family of probability measures
$\mu(p,\rho'',r'',v,w)$ supported on configurations of balls of
radius $r'':=vr'$ and having density $\rho'':=\rho'/v$. Applying
again Lemma~\ref{l:proc-relations}(4) we see that for each
$\rho'',r',v,w$ the measure $\mu(p,\rho'',r'',v,w)$ is
$\pi(p,v,r'',\IR)$-invariant.

On the other hand, these measures cannot be massive, since they are
supported by very thin sets. To construct a massive (albeit
non-ergodic) invariant measure from this family we consider a
measure $\mu(p,\rho'',r'',v)$ having marginals
$\mu(p,\rho'',r'',v,w)$ on sub-lattices $v\IZ+w$ and uniformly
distributed with respect to the parameter $w\in[0,v)$.

The construction of the massive invariant measure for
$\pi(p,v,r,\IZ)$ with general $v\in\IZ_+$ is exactly the same except 
that $w\in\{0,1,\dots,v-1\}$ and $r"=1/2$ or $r"=0$.

\subsection{Stochastic stability}\label{s:stab}

Despite very substantial differences in the constructions of 
nontrivial invariant measures in deterministic and random 
situations the structure of the invariant measure in both 
cases is described in terms of $2\times2$ Markov matrices 
$\{p_{ij}^{(p)}\}$. Therefore to prove the stochastic stability 
one only needs to check that 
$$p_{ij}^{(p)}\toas{p\to1}p_{ij}^{(1)}.$$ 
To this end one uses explicit formulas for the entries of 
the Markov matrices. Namely by (\ref{e:p10})
$$ p_{10}^{(p)}=(1-a)/(1-pa) \toas{p\to1} 1 ,$$
which coincides with the corresponding entry in the deterministic case, 
while $p_{11}^{(p)} \toas{p\to1} 0$. 
Passing to the limit as $p\to1$ one gets the corresponding limit 
relations for two other entries $p_{00},p_{01}$, which define 
uniquely the Markov measure in the lattice deterministic setting 
with $v=1$. Thus the nontrivial invariant measures for the process 
$\pi(p=1,v=1,r=1/2,\IZ)$ are stochastically stable.  

Recall now that the nontrivial invariant measures for the exclusion
 type processes in continuum $\pi(p,v,r,\IR)$ were constructed in 
Sections~\ref{s:massiv-det},\ref{s:massiv} in three steps 
through changes of variables which do not depend on the choice 
of the variable $p$. Therefore this construction withstand the limit 
transition as $p\to1$, which proves the stochastic stability for 
the general process $\pi(p=1,v,r,\IR)$. Similarly one proves 
the stochastic stability for the lattice processes with long jumps 
$\pi(p=1,v\in\IZ_+,r\in\{0,1/2\},\IZ)$.

Let me note that the Markov measures are unique among 
stochastically stable translationally invariant measures, but 
there are non-translationally invariant ones. To demonstrate 
this consider a $\delta$-measure $\mu$ supported by a 
single configuration $x:=(\dots,0,x_0=0,x_1=1,1,\dots)$. 
Then this configuration is a fixed point for the process 
$\pi(p,v=1,r=1/2,\IZ)$ for each $0<p\le1$ and hence the 
measure $\mu$ is stochastically stable.

\section{Average velocities}\label{s:vel}

Explicit expression for the invariant measure $\mu$ allows one 
to derive by Birkhoff's ergodic theorem the formula for the average 
velocity for $\mu$-a.a. initial configurations. In what follows we are 
interested in a much more broad setting of all particle configurations 
having well defined densities. This extension can be justified as 
follows. Theorem~\ref{t:vel-part} shows that for all particle configurations 
of a given density the average velocity is the same (or for all of these 
configurations the average velocity is not well defined). 
Therefore it is enough to calculate this statistics for a most suitable 
single initial configuration of given density. Choosing a configuration 
typical with respect to the invariant measure $\mu$ we achieve this goal. 
Since this construction does not depend on specific properties of the 
invariant measure we shall not repeat this argument in further calculations.

\subsection{Lattice TASEP with $v=1$}\label{s:erg-1}

Let us use the constructed Markov measures to calculate the average
velocities for the process $\pi(p,v,r=1/2,\IZ)$.

\begin{lemma}\label{l:vel-1}
$V(\rho)=(1-\sqrt{1-4p\rho(1-\rho)})/(2\rho)$ for $\pi(p<1,v,r=1/2,\IZ)$.
\end{lemma}

\proof For a given configuration of particles, a particle may move iff
the next site of the lattice is not occupied, i.e. only in the situation $10$.
Therefore the average velocity is equal to the stationary probability of
the jump to the right, which in turn is equal to $pp_{10}$. The representation 
of the Markov invariant measure obtained earlier immediately gives
the formula for the average velocity %
$$V(\rho,p,1,1/2)=pp_{10}=p(1-a)/(1-pa)\in[0,p] .$$ 
Substituting the value
$a=a(\rho)$ according to the formula (\ref{e:par}) we get %
\beq{e:vel0}{V(\rho,p,v=1,r=1/2)
             =\left[1-\sqrt{1-4p\rho(1-\rho)}\right]/(2\rho).} %
\qed

\subsection{Process  in continuum}\label{s:erg-3}

By means of results of Theorem~\ref{t:vel-part}  and
Lemma~\ref{l:proc-relations} the relation~(\ref{e:vel0})
can be transfered to the processes of the 3d type $\pi(p,v,r,\IR)$ 
without changes.
Note that the formula~(\ref{e:vel0}) was already known in physics
publications for the lattice processes of type 1 with $v=1$ (see \cite{SS}).

It is important to say that the naive transition from $v=1$ to
$v>1$ directly in the class of lattice processes (using the
invariance of sub-lattices with step multiple to $v$) is impossible
(or rather so we can study only low-density $<1/v$ configurations).
Instead, we use the self-similarity of processes of type 3 acting in 
a continuous space.

\begin{lemma}\label{l:sim}
Let the $\pi(p,v,r,\IR)$ type processes $x^t,\2x^t$ with parameters
$r=\2r=0,~v>0,~\2v=uv,~p=\2p\in(0,1)$ having initial configurations
$\2x^0=ux^0$ with some $u>0$ be statically coupled.
Assume also that $\rho(x)$ and $V(x)$ be well defined. Then %
\beq{e:den-vel}{\rho(\2x)=\rho(x)/u,\quad \2V(\2x)=uV(x).}
\end{lemma}
\proof By a straightforward calculation. \qed

Applying these similarity transformations to the special case
described in the relation (\ref{e:vel0}) we obtain the general formula
for the average velocity (\ref{e:vel}).

\begin{corollary}\label{c:1-3}
\bea{ V(\rho,p,v>0,r=0) \a=v V(v\rho,p,v=1,r=1/2) \\
 \a= \frac{1+v\rho-\sqrt{(1+v\rho)^2-4pv\rho}}{2\rho}\\
 \a\toas{p\to1}\min(v,~1/\rho). }
\end{corollary}

Since the last term above coincides with the average velocity for the 
deterministic traffic map (see \cite{Bl10}), the limit transition in the last 
relation shows that the average velocity is stochastically stable. 

Indeed, under a spatial change of variables $1\to v$ and the transition
from the configuration of balls of radius $r=1/2$ to the configuration
of point-particles (i.e. $\2r=0$) with the same sequence of gaps we get
$$\2\rho=\rho/(1-2r\rho)=\rho/(1-\rho),$$ 
hence $\rho=\2\rho/(1+\2\rho)$ and
\bea{ V(\2\rho,p,v=1,\2r=0) \a=\frac{1-\sqrt{1-4p\rho(1-\rho)}}{2\rho} \\
 \a= \frac{1-\sqrt{1-4p\frac{\2\rho}{1+\2\rho}(1-\frac{\2\rho}{1+\2\rho})}}
           {2\frac{\2\rho}{1+\2\rho}} \\
 \a= \frac{1+\2\rho-\sqrt{(1+\2\rho)^2-4p\2\rho}}{2\2\rho} .}
Therefore by Lemma~\ref{l:sim}
\bea{ V(\rho,p,v>0,r=0) \a=v V(v\rho,p,v=1,r=1/2) \\
 \a= \frac1v\times \frac{1+v\rho-\sqrt{(1+v\rho)^2-4pv\rho}}{2v\rho} \\
 \a= \frac{1+v\rho-\sqrt{(1+v\rho)^2-4pv\rho}}{2\rho}
 \toas{p\to1}\min(v,~1/\rho) .}

This finishes the proof of Theorem~\ref{t:vel-main}.

\subsection{Lattice exclusion processes  with long jumps}\label{s:erg-1-2}

The results for $\pi(p,v,r,\IR)$ with arbitrary $v\in\IZ_+$ are transferred
directly by Lemma~\ref{l:proc-relations} back to the lattice cases.
Indeed, Theorem~\ref{t:vel-part} shows that it is enough to derive
the formula for the average velocity for a specially chosen initial
configuration of a given density. On the other hand, by
Lemma~\ref{l:proc-relations}(1 and 2)  a realization of the
processes $\pi(p,v,r=1/2,\IZ)$ and $\pi(p,v,r=0,\IZ)$ starting
from certain configurations coincide with a realization of the
$\pi(p,v,r,\IR)$ process statically coupled to the lattice process and
starting from the same initial configuration.

It is worth noting that a naive application of the property described 
in Lemma~\ref{l:proc-relations}(4) seems to extend the results
about the process $\pi(p,v=1,r=1/2,\IZ)$ directly to
$\pi(p,v>1,r=1/2,\IZ)$ restricting the latter process to invariant
sub-lattices $v\IZ+w,~w=0,1,\dots,v-1$. A close look shows that
this is indeed the case but only for configurations of low density
$\rho<1/v$, since otherwise particles from the same configuration 
located at different sub-lattices will interact.

\subsection{Heterogeneous particles (of different sizes)}\label{s:het-size}

Thinking about the processes under consideration as models of traffic
flows it is reasonable to take into account that vehicles need not to
be of the same size. From this point of view we consider an exclusion
type process in continuum with particle configurations consisting of balls
with varying sizes, i.e. the radius of the $i$-th ball is equal to $r_i\ge0$.
Our aim is to show that the dependence of the average velocity on density
in this case can be easily obtained from the corresponding result for
the case of balls of the same radius.

Let $x$ be a bi-infinite admissible configuration of particles represented
by balls of radiuses $r_i\ge0$ with the average value
$$\b{r}:=\lim_{n\to\infty}\frac1n\sum_{i=0}^{n-1}r_i$$ and centered
at points $x_i\in\IR$.
The notion of admissibility and the law of dynamics should be slightly
rewritten (in comparison to the homogeneous case):
$$x_i + r_i\le x_{i +1}-r_{i+1} ,$$
$${x_i^{t+1}=\function{\min\{x_i^t+v,x_{i+1}^t-r_i-r_{i+1}\}
                          &\mbox{with probability } p \\
                     x_i^t &\mbox{with probability } 1-p} .}$$

\begin{theorem}\label{t:heter-sizes} Let the process $x^t$ defined
above and the $\pi(p,v,\2r,\IR)$ process $\2x^t$ be statically
coupled and let $\2r=\b{r}$. Then the average velocities of these
processes coincide.
\end{theorem}
\proof A simple generalization of the affine conjugation between
configurations of different ball's sizes introduced in
Lemma~\ref{l:proc-relations} allows one to make the bijection 
between the configurations $x^t,\2x^t$.
Indeed, consider an affine map defined by the relation
$$ (\phi(\2x))_i := x_i - 2\sum_{j=0}^i(r_j-\2r) .$$
Then
\bea{ \frac1n(x_{n-1}^t-x_0^t)
 \a= \frac1n( (\phi^{-1}(\2x^t))_{n-1} -  (\phi(\2x^t))_0) \\
 \a= \frac1n(\2x_{n-1}^t-\2x_0^t)
    + \frac2n\sum_{j=0}^{n-1}(r_j-\2r) - \frac2n(r_0-\2r) \\
\a\toas{n\to\infty}1/\rho(\2x^t) + 2(\b{r}-\2r) .}
Thus $\rho(x^t)=\rho(\2x^t)$ if $\2r=\b{r}$.

Now we are ready to calculate the average velocity $V(x)$:
\bea{ V(x,i,t) \a= \frac1t(x_i^t - x_i^0)
 =  \frac1t( (\phi^{-1}(\2x^t))_i -  (\phi(\2x^0))_i) \\
 \a=  \frac1t(\2x_i^t - \2x_i^0) + \frac2t\sum_{j=0}^i(r_j-\2r) - \frac2t(r_0-\2r) \\
\a\toas{t\to\infty}V(\2x) + 2(\b{r}-\2r) ,}
which proves our claim. \qed

\subsection{Heterogeneous space}\label{s:het-space}

So far we have considered only exclusion processes acting on
homogeneous spaces. In \cite{Bl11} we introduced and studied
a modification of the deterministic version of the exclusion process
in continuum which takes into account the presence of static
obstacles (traffic lights) for the motion of point particles
(i.e. $r=0$). Fix an arbitrary point-particle configuration
$z=(z_j)_{j\in\IZ}\in X(0,\IR)$ whose elements correspond to
positions of obstacles. Then the formula (\ref{e:dyn}) can be
rewritten as follows: %
\beq{e:dyn1}{x_i^{t+1}=\function{\min\{x_i^t+v,x_{i+1}^t,z_{j(x_i^t)}\}
                          &\mbox{with probability } p \\
                     x_i^t &\mbox{with probability } 1-p} ,}%
where $j(x_i^t):=\min\{k\in\IZ:~x_i^t\le z_k\}$. Thus the
``obstacles'' suspend the movement of particles, taking into account
the time necessary to overtake an obstacle.

For a given $v>0$ and a configuration of obstacles $z$ denote by
$\t{z}$ the {\em extended} configuration of obstacles obtained by
inserting between each pair of entries $z_i,z_{i+1}$ new
$\intp{(z_{i+1}-z_i)/v}$ `virtual' obstacles at distances
$v$ between them starting from the point $z_i$. Here
$\intp{u}$ stands for the integer part of the number $u$.

\begin{theorem} \label{t:vel-obs}
For given $v>0<p<1$ and any configurations $x,z\in X$ for which
the densities $\rho(x), \rho(\t{z})$ are well defined %
\beaq{e:vel-obs0} {V(x,z)\a=\frac{\rho(x)+\rho(\t{z})
        -\sqrt{(\rho(x)+\rho(\t{z}))^2-4p\rho(x)\rho(\t{z})}}
        {2\rho(x)\rho(\t{z})}  \cr
        \a\toas{p\to1}\min\{1/\rho(\t{z}),~1/\rho(x)\}. }%
\end{theorem}

Note that the average velocity in the above formula does not
depend explicitly on the local velocity $v$, however the latter is
included to the construction of the extended configuration
$\t{z}$, in particular $\rho(\t{z})\geq1/v$.
It is interesting to note also that in \cite{Bl11} it was shown that under
more general setting with the non-degenerate distribution of random
iid local velocities $v_i^t$ the average particle velocity may not exist.

The difficulty of the analysis here is is that the inhomogeneity of
the space in which the collective random walk takes place
(the presence of obstacles), generally does not permit the existence
of invariant measures.\footnote{For the existence of invariant
     measures one needs at least the condition of stationarity
     for the configurations of obstacles $z$.} %
Therefore, the main step of the approach we used -- the
construction of a massive invariant measure is impossible.
The complete proof of this result needs a modification of
the dynamical coupling construction elaborated in \cite{Bl10,Bl11}
and will be discussed elsewhere. Here we only give an idea of the proof.

Using the technique developed in \cite{Bl11} for the deterministic
version of the problem it can be shown that the calculation of
average velocities $V(x,z,v,p)$ can be reduced to the analysis of
a Markov process of type~2, acting (in contrast to the
already-studied setting) on {\em a inhomogeneous lattice}
$\bR:=\tilde z$. The existence of the density of the configuration
$\tilde z$ allows one to transfer the results obtained for the conventional
integer lattice $\IZ$ to the inhomogeneous case under consideration.
This completes the construction.

\newpage
{

}
\end{document}